\documentclass{article}
\usepackage{amsfonts}
\usepackage{makeidx}
\usepackage{amssymb}
\usepackage{amsmath}
\usepackage{amstext}
\usepackage{graphicx}

\newtheorem{theorem}{Theorem}

\newtheorem{corollary}[theorem]{Corollary}

\newtheorem{definition}[theorem]{Definition}

\newtheorem{lemma}[theorem]{Lemma}

\newtheorem{proposition}[theorem]{Proposition}
\newtheorem{remark}[theorem]{Remark}

\input tcilatex
\begin{document}

\title{A non-archimedean algebra and the Schwartz impossibility theorem}
\author{Vieri Benci\thanks{%
Department of Mathematics, University of Pisa, Via F. Buonarroti 1/c, 56127
Pisa, ITALY and Department of Mathematics, College of Science, King Saud
University, Riyadh, 11451, SAUDI ARABIA. e-mail: \texttt{benci@dma.unipi.it}}
\and Lorenzo Luperi Baglini\thanks{%
University of Vienna, Faculty of Mathematics, Oskar-Morgenstern-Platz 1,
1090 Vienna, AUSTRIA, e-mail: \texttt{lorenzo.luperi.baglini@univie.ac.at},
supported by grant P25311-N25 of the Austrian Science Fund FWF.}}
\maketitle

\begin{abstract}
In the 1950s L. Schwartz proved his famous impossibility result: for every $%
k\in \mathbb{N}$ there does not exist a differential algebra $(\mathfrak{A}%
,+,\otimes ,D)$ in which the distributions can be embedded, where $D$ is a
linear operator that extends the distributional derivative and satisfies the
Leibnitz rule (namely $D(u\otimes v)=Du\otimes v+u\otimes Dv$) and $\otimes $
is an extension of the pointwise product on $\mathcal{C}^{0}(\mathbb{R})$.

In this paper we prove that, by changing the requests, it is possible to
avoid the impossibility result of Schwartz. Namely we prove that it is
possible to construct an algebra of functions $(\mathfrak{A,}+,\otimes ,D)$
such that (1) the distributions can be embedded in $\mathfrak{A}$ in such a
way that the restriction of the product to $\mathcal{C}^{1}(\mathbb{R})$
functions agrees with the pointwise product, namely for every $f,g\in 
\mathcal{C}^{1}(\mathbb{R})$%
\begin{equation*}
\Phi (fg)=\Phi (f)\otimes \Phi \left( g\right) ,
\end{equation*}%
and (2) there exists a linear operator $D:\mathfrak{A}\rightarrow \mathfrak{A%
}$ that extends the distributional derivative and satisfies a weak form of
the Leibnitz rule.

The algebra that we construct is an algebra of restricted ultrafunctions,
which are generalized functions defined on a subset $\Sigma $ of a
non-archimedean field $\mathbb{K}$ (with $\mathbb{R}\subset \Sigma \subset 
\mathbb{K}$) and with values in $\mathbb{K}$. To study the restricted
ultrafunctions we will use some techniques of nonstandard analysis.


\medskip

\end{abstract}

\tableofcontents

\section{Introduction}

There is an issue regarding distributions that is important for a variety of
applications, namely the problem of defining a multiplication of
distributions that satisfies some property of coherence with respect to the
weak derivative and to the restriction to continuous functions (see \cite%
{libro}, Chapter 1 for a discussion on this topic). A possible way to define
such a multiplication is to embed the space of distributions in a
differential algebra $(\mathfrak{A},+,\otimes ,D)$ and to use $\otimes $ to
define the multiplication of distributions. A famous result that limits this
approach was proved by L. Schwartz in \cite{Schwartz}: he proved that it is
impossible to construct a differential algebra $(\mathfrak{A},+,\otimes ,D)$
such that

(i) there is a linear embedding%
\begin{equation*}
\Phi :\mathcal{D}^{\prime }(\mathbb{R})\rightarrow \mathfrak{A}
\end{equation*}%
such that $\Phi \left( 1\right) $ is the unity in $\mathfrak{A};$

(ii) there is a linear operator $D:$ $\mathfrak{A}\rightarrow \mathfrak{A}$
such that the following diagram%
\begin{equation*}
\begin{array}{ccc}
\mathcal{D}^{\prime }(\mathbb{R}) & \overset{\partial }{\longrightarrow } & 
\mathcal{D}^{\prime }(\mathbb{R}) \\ 
\Phi \downarrow &  & \Phi \downarrow \\ 
\mathfrak{A} & \overset{D}{\longrightarrow } & \mathfrak{A}%
\end{array}%
\end{equation*}%
commutes, where $\partial $ is the usual distributional derivative;

(iii) the restriction of $\otimes $ to the continuous functions agrees with
the pointwise product, namely%
\begin{equation*}
\Phi (fg)=\Phi (f)\otimes \Phi \left( g\right) ;
\end{equation*}

(iv) the Leibnitz rule holds:%
\begin{equation*}
D\left( uv\right) =Duv+uDv.
\end{equation*}

Actually, for every $k\in \mathbb{N},$ the impossibility result holds even
if we modify (iii) as follows:

\bigskip

(iii)$_{k}$ the restriction of $\otimes $ to $\mathcal{C}^{k}(\mathbb{R}%
)\times \mathcal{C}^{k}(\mathbb{R})$ agrees with the pointwise product,
namely%
\begin{equation*}
\Phi (fg)=\Phi (f)\otimes \Phi \left( g\right) .
\end{equation*}

In order to embedd the distributions in a differential algebra one has to
weaken at least one of the requests (i),..., (iv). A famous approach to this
problem is given by Colombeau's Algebras, in which (iii) is replaced by

\bigskip

(iii)$_{\infty }$ the restriction of $\otimes $ to $\mathcal{C}^{\infty }(%
\mathbb{R})\times \mathcal{C}^{\infty }(\mathbb{R})$ agrees with the
pointwise product, namely%
\begin{equation*}
\Phi (fg)=\Phi (f)\otimes \Phi \left( g\right) .
\end{equation*}

Jean-Fran\c{c}ois Colombeau proved the existence of algebras satisfying (i),
(ii), (iii)$_{\infty }$, (iv). The central ideas of his construction were
first published in \cite{Col82}, \cite{Col83a} and \cite{Col84}, and the
foundations of his work are written in the books \cite{Col84}, \cite{Col85}.
For a more recent reference on this topic we suggest the book \cite{libro}.

In this paper we prove a different existence result by relaxing the requests
(i), (ii), (iii), (iv) in a different way. We slightly weaken (iii) but we
weaken (iv) in a more substantial way. We substitute (iii) with (iii)$_{1}$,
namely

\bigskip

(iii)$_{1}$: the restriction of $\otimes $ to $\mathcal{C}^{1}(\mathbb{R}%
)\times \mathcal{C}^{1}(\mathbb{R})$ agrees with the pointwise product,
namely%
\begin{equation*}
\Phi (fg)=\Phi (f)\otimes \Phi \left( g\right) .
\end{equation*}

\bigskip

Let us show how we weaken the Leibnitz rule (iv). If $u$ and $v$ are
functions. by integrating (iv) we get 
\begin{equation}
\int Duv+\int uDv=\left[ uv\right] _{-\infty }^{+\infty },  \label{lisa}
\end{equation}%
provided that%
\begin{equation*}
\left[ uv\right] _{-\infty }^{+\infty }=\ \underset{x\rightarrow +\infty }{%
\lim }u(x)v(x)-\ \underset{x\rightarrow -\infty }{\lim }u(x)v(x).
\end{equation*}%
is well defined. Clearly (\ref{lisa}) is a weaker request than the Leibnitz
rule. We make a request on the elements of $\mathfrak{A}$ which generalizes (%
\ref{lisa}):

\bigskip

(iv)' (\textbf{Weak Leibnitz Rule}) For every $u,v\in \mathfrak{A}$%
\begin{equation}
\left\langle Du,v\right\rangle +\left\langle u,Dv\right\rangle \ =\left[ uv%
\right] _{-\beta }^{+\beta },  \label{isa}
\end{equation}%
where $\left\langle u,v\right\rangle $ is a scalar product such that, for
every $f,g\in \mathcal{C}_{0}^{1}(\mathbb{R}),$%
\begin{equation*}
\left\langle \Phi (T_{f}),\Phi (T_{g})\right\rangle =\int f(x)g(x)dx
\end{equation*}%
where $T_{f}$, $T_{g}$ are the distributions associated to $f,g$ and $\beta $
is a suitable "point at infinity".

\bigskip

Notice that (\ref{lisa}) is used to define the notion of weak derivative and
the duality in the theory of distribution. So, even if (\ref{lisa}) and (\ref%
{isa}) are weaker than the Leibnitz rule, they are essential in the
applications.

We will show that the requests (i), (ii), (iii)$_{1}$, (iv)' are consistent
by constructing explicitly an algebra $\mathfrak{A}$ that satisfies these
properties. This construction will be done by using the space of
ultrafunctions, which is a space of generalized functions that has been
introduced in \cite{ultra} and further studied in \cite{belu2012}\ and \cite%
{belu2013}. An interesting feature of the algebra $\mathfrak{A}$ is that
there exists a non-archimedean field $\mathbb{K\supset R}$ such that $%
\mathfrak{A}$ is a subalgebra of the algebra of functions%
\begin{equation*}
u:\Sigma \rightarrow \mathbb{K\ \ }\text{where}\mathbb{\ \ R}\subset \Sigma
\subset \mathbb{K},
\end{equation*}%
equipped with the pointwise operations:%
\begin{equation*}
(u+v)(x)=u(x)+v(x);\ \ (u\otimes v)(x)=u(x)v(x).
\end{equation*}

Our construction uses some tools of nonstandard analysis. In the literature,
nonstandard analysis has been used many times to study questions related to
Schwartz's impossibility result and to the Colombeau's algebras. For
example, in \cite{Rob2}, the field of asymptotic real numbers has been
introduced, which is related to Colombeau algebras; also we recall the more
recent results in \cite{Obe} and \cite{Tod}. However, our construction is
quite different with respect to these previous nonstandard approaches.

\subsection{Notations and definitions\label{not}}

We use this section to fix some notations and to recall some definitions:

\begin{itemize}
\item $\mathfrak{F}\left( X,Y\right) $ denotes the set of all functions from 
$X$ to $Y;$

\item $\mathfrak{F}\left( \mathbb{R}\right) =\mathfrak{F}\left( \mathbb{R},%
\mathbb{R}\right) ;$

\item $\mathcal{C}\left( \mathbb{R}\right) $ denotes the set of continuous $%
f:\mathbb{R\rightarrow R};$

\item $\mathcal{C}_{0}\left( \mathbb{R}\right) $ denotes the set of
functions in $\mathcal{C}\left( \mathbb{R}\right) $ having compact support;

\item $\mathcal{C}^{k}\left( \mathbb{R}\right) $ denotes the set of
functions in $\mathcal{C}\left( \mathbb{R}\right) $ which have continuous
derivatives up to the order $k;$

\item $\mathcal{C}_{0}^{k}\left( \mathbb{R}\right) $ denotes the set of
functions in $\mathcal{C}^{k}\left( \mathbb{R}\right) \ $having compact
support;

\item $\mathcal{D}\left( \mathbb{R}\right) $ denotes the set of the
infinitely differentiable functions with compact support$;\ \mathcal{D}%
^{\prime }\left( \mathbb{R}\right) $ denotes the topological dual of $%
\mathcal{D}\left( \mathbb{R}\right) $, namely the set of distributions on $%
\mathbb{R};$

\item if $\mathbb{K\ }$is a linearly ordered field and $a,b\in \mathbb{K},$
then

\begin{itemize}
\item $\left[ a,b\right] _{\mathbb{K}}=\{x\in \mathbb{K}:a\leq x\leq b\};$

\item $\left( a,b\right) _{\mathbb{K}}=\{x\in \mathbb{K}:a<x<b\};$
\end{itemize}

\item an element $k$ of an ordered field $\mathbb{K}$ is infinite if $|k|>n$
for every natural number $n;$

\item an ordered field $\mathbb{K}$ is non-archimedean if it contains
infinite elements;

\item a field $\mathbb{K}$ is superreal if it properly contains the field $%
\mathbb{R}$.
\end{itemize}

\section{\textbf{The main result}}

In this section we state the main result of the paper, which will be proved
in section \ref{dimostraz}.

\begin{theorem}
\label{gonerilla}There exists an algebra $(\mathfrak{A},+,\cdot ,\mathrm{D})$
that satisfies the following properties:

\begin{itemize}
\item ($\mathfrak{A}$-0) (\textbf{Algebraic structure}) $\mathfrak{A}%
\subseteq \mathfrak{F}\left( \Sigma ,\mathbb{K}\right) $ where $\mathbb{K}$
is a non-archimedean field and $\Sigma $ is a set such that%
\begin{equation*}
\mathbb{R}\subset \Sigma \subset \mathbb{K};
\end{equation*}%
$\mathfrak{A}$ is an algebra equipped with the pointwise operations:%
\begin{equation*}
(u+v)(x)=u(x)+v(x);\ \ (u\cdot v)(x)=u(x)\cdot v(x).
\end{equation*}

\item ($\mathfrak{A}$-1) (\textbf{Embedding of distributions}) There is a
linear embedding%
\begin{equation*}
\Phi :\mathcal{D}^{\prime }(\mathbb{R)}\rightarrow \mathfrak{A}
\end{equation*}%
and a bilinear form $\left\langle \cdot ,\cdot \right\rangle :\mathfrak{A}%
\times \mathfrak{A}\rightarrow \mathbb{K}$ such that, $\forall T\in \mathcal{%
D}^{\prime }(\mathbb{R)}$, $\forall \varphi \in \mathcal{D(}\mathbb{R)}$, 
\begin{equation*}
T[\varphi ]=\left\langle \Phi \left( T\right) ,\Phi (T_{\varphi
})\right\rangle .
\end{equation*}

\item ($\mathfrak{A}$-2) (\textbf{Extension of the derivative}%
\QTR{frametitle}{) }There is a linear operator $\mathrm{D}:$ $\mathfrak{A}%
\rightarrow \mathfrak{A}$ such that the diagram%
\begin{equation}
\begin{array}{ccc}
\mathcal{D}^{\prime }(\mathbb{R)} & \overset{\partial }{\longrightarrow } & 
\mathcal{D}^{\prime }(\mathbb{R)} \\ 
\downarrow \Phi &  & \downarrow \Phi \\ 
\mathfrak{A} & \overset{\mathrm{D}}{\longrightarrow } & \mathfrak{A}%
\end{array}
\label{bello}
\end{equation}%
commutes, where $\partial $ is the usual distributional derivative.

\item ($\mathfrak{A}$-3) (\textbf{Extension of the product}) The restriction
of $\cdot $ to $\mathcal{C}^{1}(\mathbb{R)}$ agrees with the pointwise
product namely, if $f,g\in \mathcal{C}^{1}(\mathbb{R)},\ $then%
\begin{equation*}
\Phi (T_{fg})=\Phi (T_{f})\cdot \Phi \left( T_{g}\right) .
\end{equation*}

\item ($\mathfrak{A}$-4) (\textbf{Weak Leibnitz rule}) For every $u,v\in 
\mathfrak{A}$ the following holds:%
\begin{equation*}
\left\langle \mathrm{D}u,v\right\rangle +\left\langle u,\mathrm{D}%
v\right\rangle \ =\left[ uv\right] _{-\beta }^{+\beta },
\end{equation*}%
where $\beta =\max (\Sigma ),$ $-\beta =\min (\Sigma ).$

\item ($\mathfrak{A}$-5) (\textbf{Locality of the extension}\QTR{frametitle}{%
) }If the support of a distribution $T\ $is included in $\left[ a,b\right]
\subset \mathbb{R}$ then, for every $x\in \Sigma \backslash \left[ a,b\right]
_{\mathbb{K}},$ we have 
\begin{equation*}
\Phi (T)(x)=0.
\end{equation*}
\end{itemize}
\end{theorem}

Let us observe that, since $\beta \in \mathbb{K\setminus R}$, $\beta $ is an
infinite number in $\mathbb{K}$\ and that every algebra given by Theorem \ref%
{gonerilla} satisfies the requests (i), (ii), (iii)$_{1}$, (iv)' outlined in
the introduction; moreover, as an immediate consequence of Theorem \ref%
{gonerilla}, the operator $\mathrm{D}$ and the scalar product $\left\langle
\cdot ,\cdot \right\rangle $ have properties similar to the duality of
distributions. In the following corollary we identify every function $f\in $ 
$\mathcal{C}_{0}^{1}(\mathbb{R)}$ with its counterpart in $\mathfrak{A}$,
namely with $\Phi (T_{f}).$

\begin{corollary}
\label{firuli}$\forall u\in \mathfrak{A,}$ $\forall f\in $ $\mathcal{C}%
_{0}^{1}(\mathbb{R)},$ $\left\langle \mathrm{D}u,f\right\rangle
=-\left\langle u,\partial f\right\rangle .$
\end{corollary}

\textbf{Proof. }By ($\mathfrak{A}$-4) we have $\left\langle \mathrm{D}%
u,f\right\rangle +\left\langle u,\mathrm{D}f\right\rangle \ =\left[ uf\right]
_{-\beta }^{+\beta },$ and by ($\mathfrak{A}$-5) we get that $[uf]_{-\beta
}^{\beta }=0.$ Moreover, by ($\mathfrak{A}$-2) it follows that $\mathrm{D}%
f=\Phi (\partial T_{f})=\partial f$ (with respect to our identification),
hence we can conclude. $\square $

\section{Construction of the Algebra}

\subsection{The Ultrafunctions}

Throughout this section we assume that the reader has a basic knowledge of
nonstandard analysis (for a general reference on the subject, see e.g. \cite%
{rob}). We work in a (at least) $\left( 2^{\mathfrak{c}}\right) ^{+}$%
-saturated extension of the real numbers\footnote{%
We recall that, given a cardinal number $\mathit{k}$, a nonstandard model
has the $\mathit{k}^{+}$-saturation property if for every family $\mathfrak{F%
}$ of internal sets with the finite intersection property and with $|%
\mathfrak{F}|\leq \mathit{k}$ the intersection $\bigcap_{A\in \mathfrak{F}}A$
is not empty.} (where $\mathfrak{c}$ stands for the cardinality of
continuum), and we take as standard universe the superstructure $V(\mathbb{R}%
)$ on $\mathbb{R}$. We recall that, given a set $A$ in $V(\mathbb{R}),$ $%
A^{\sigma }$ is the set%
\begin{equation*}
A^{\sigma }=\{a^{\ast }\mid a\in A\}.
\end{equation*}%
We let $\Lambda $ denote a hyperfinite set in $\mathfrak{F}(\mathbb{R}$,$%
\mathbb{R})^{\ast }$ with $\mathfrak{F}(\mathbb{R}$,$\mathbb{R})^{\sigma
}\subseteq \Lambda $. We let%
\begin{equation*}
\widetilde{\mathcal{C}^{1}}(\mathbb{R})=Span\{\mathcal{C}^{1}(\mathbb{R}%
)^{\ast }\cap \Lambda \}.
\end{equation*}%
Let us observe that, by definition, $\widetilde{\mathcal{C}^{1}}(\mathbb{R})$
is an internal vector space of hyperfinite dimension and $\mathcal{C}^{1}(%
\mathbb{R})^{\sigma }\subseteq \widetilde{\mathcal{C}^{1}}(\mathbb{R}).$

\begin{definition}
Let $\beta $ be a positive infinite number. We call ultrafunctions the
elements of the space $V_{\Lambda },$ where%
\begin{equation*}
V_{\Lambda }=\{u_{\upharpoonleft _{\lbrack -\beta ,\beta ]}}\mid u\in 
\widetilde{\mathcal{C}^{1}}(\mathbb{R})\}.
\end{equation*}
\end{definition}

\begin{remark}
In our previous works (\cite{belu2012}, \cite{belu2013}) we called $%
\widetilde{\mathcal{C}^{1}}(\mathbb{R})$ the space of ultrafunctions
generated by $\mathcal{C}^{1}(\mathbb{R})$ (which was constructed in a
different, but equivalent, way). In this paper we slightly changed our
definition of "ultrafunction".
\end{remark}

From now on, with some abuse of notation, we will say that a function $%
\varphi $ is in $V_{\Lambda }$ meaning that the restriction $\varphi
_{\upharpoonleft _{\lbrack -\beta ,\beta ]}}\in V_{\Lambda }.$ Similarly,
when we say that $f^{\ast }\in V_{\Lambda }$ we mean that $f_{\upharpoonleft
_{\lbrack -\beta ,\beta ]}}^{\ast }\in V_{\Lambda }$.

On the space $V_{\Lambda }$ we can define a notion of derivative by duality
as follows:

\begin{definition}
For every ultrafunction $u\in V_{\Lambda },$ the derivative $Du$ of $u$ is
the unique ultrafunction such that, for every $v\in V_{\Lambda },$%
\begin{equation*}
\int_{-\beta }^{\beta }Du(x)v(x)dx=\int_{-\beta }^{\beta }\partial ^{\ast
}u(x)v(x)dx,
\end{equation*}%
where $\int_{-\beta }^{\beta }$ denotes the extension of the Lebesgue
integral to $\mathbb{R}^{\ast }$ with limits $-\beta ,\beta .$
\end{definition}

Let $P_{V_{\Lambda }}:\mathcal{C}^{0}(\mathbb{R})^{\ast }\rightarrow
V_{\Lambda }$ be the orthogonal projection w.r.t. the $L^{2}$ scalar product
defined on $[-\beta ,\beta ]$, namely for every $f\in \mathcal{C}^{0}(%
\mathbb{R})^{\ast }$ $P_{V_{\Lambda }}f$ is the unique ultrafunction such
that, for every ultrafunction $u$, we have 
\begin{equation*}
\int_{-\beta }^{\beta }f(x)u(x)dx=\int_{-\beta }^{\beta }P_{V_{\Lambda
}}f(x)u(x)dx.
\end{equation*}%
Then $D$ can be equivalently expressed by composition as follows:%
\begin{equation*}
D=P_{V_{\Lambda }}\circ \partial ^{\ast }.
\end{equation*}%
An immediate consequence of the definition is that, if $f\in \mathcal{C}^{2}(%
\mathbb{R}),$ then 
\begin{equation*}
Df^{\ast }=(\partial f)^{\ast }.
\end{equation*}%
In fact, if $f\in \mathcal{C}^{2}(\mathbb{R})$ then $\partial f\in \mathcal{C%
}^{1}(\mathbb{R})$ and, since $\mathcal{C}^{1}(\mathbb{R})^{\sigma
}\subseteq $ $\widetilde{\mathcal{C}^{1}}(\mathbb{R}),$ we have $(\partial
f)^{\ast }=P_{V_{\Lambda }}(\partial f)^{\ast }=Df^{\ast }.$

For our aims, the most important property of $D$ is the following:

\begin{theorem}
\label{facundo}For every $u,v\in V_{\Lambda }$ we have 
\begin{equation*}
\int_{-\beta }^{\beta }Du(x)v(x)dx=-\int_{-\beta }^{\beta
}u(x)Dv(x)dx+[uv]_{-\beta }^{\beta }.
\end{equation*}
\end{theorem}

\textbf{Proof. }Let us compute $\int_{-\beta }^{\beta }Du(x)v(x)dx:$%
\begin{eqnarray*}
\int_{-\beta }^{\beta }Du(x)v(x)dx &=&\int_{-\beta }^{\beta }\partial
u(x)v(x)dx= \\
-\int_{-\beta }^{\beta }u(x)\partial v(x)dx+[uv]_{-\beta }^{\beta }
&=&-\int_{-\beta }^{\beta }u(x)Dv(x)dx+[uv]_{-\beta }^{\beta }.\text{ \ \ \
\ }\square
\end{eqnarray*}

\bigskip

This derivative will play a central role in the construction of the algebra $%
\mathfrak{A}$. One of its important properties is presented in the following:

\begin{proposition}
\label{derivatona}For every $k\in \mathbb{N}^{\ast },$ for every $u\in
V_{\Lambda },$ for every $\varphi \in \mathcal{D}(\mathbb{R)}$ we have the
following: 
\begin{equation*}
\int_{-\beta }^{\beta }D^{k}u(x)\cdot \varphi ^{\ast
}(x)dx=(-1)^{k}\int_{-\beta }^{\beta }u(x)\partial ^{k}\varphi ^{\ast }(x)dx.
\end{equation*}
\end{proposition}

\textbf{Proof. }By internal induction on $k$: if $k=0$ there is nothing to
prove. Let us suppose the statement true for $k$. Then%
\begin{eqnarray*}
\int_{-\beta }^{\beta }D^{k+1}(u(x))\varphi ^{\ast }(x)dx &=&\int_{-\beta
}^{\beta }D(D^{k}(u(x)))\varphi ^{\ast }(x)dx= \\
&&-\int_{-\beta }^{\beta }D^{k}(u(x))D\varphi ^{\ast }(x)dx+\left[
D^{k}u\cdot \varphi ^{\ast }\right] _{-\beta }^{\beta }.
\end{eqnarray*}

Since $\varphi \in \mathcal{D}(\mathbb{R)}$ we have $\left[ D^{k}u\cdot
\varphi ^{\ast }\right] _{-\beta }^{\beta }=0.$ Moreover $D\varphi ^{\ast }=$
$\partial \varphi ^{\ast }\in \mathcal{D}(\mathbb{R)}$. So by inductive
hypothesis we have%
\begin{eqnarray*}
-\int_{-\beta }^{\beta }D^{k}(u(x))D\varphi ^{\ast }(x)dx &=&-\int_{-\beta
}^{\beta }D^{k}(u(x))\partial \varphi ^{\ast }(x)dx= \\
&&(-1)^{k+1}\int_{-\beta }^{\beta }u(x)\partial ^{k+1}\varphi ^{\ast }(x)dx,
\end{eqnarray*}%
and the thesis is proved. $\square $

\bigskip As stated in Theorem \ref{gonerilla}, we want the algebra $%
\mathfrak{A}$ to be a subalgebra of $\mathfrak{F}(\Sigma ,\mathbb{K}),$
where $\Sigma \subseteq \mathbb{K}$ and $\mathbb{K}$ is a non-archimedean
field. We fix $\mathbb{K=R}^{\ast }$, and to choose $\Sigma $ we use the
notion of "independent set of points" (which has been introduced in \cite%
{belu2013}):

\begin{definition}
\label{dede}Given a number $q\in \Omega ^{\ast },$ we denote by $\mathbf{%
\delta }_{q}(x)$ an ultrafunction in $V_{\Lambda }$ such that 
\begin{equation}
\forall v\in V_{\Lambda },\ \int_{-\beta }^{\beta }v(x)\mathbf{\delta }%
_{q}(x)dx=v(q).  \label{deltafunction}
\end{equation}%
$\mathbf{\delta }_{q}(x)$ is called Delta (or Dirac) ultrafunction centered
in $q$. \newline
A Delta-basis $\left\{ \mathbf{\delta }_{a}(x)\right\} _{a\in \Sigma }$ $%
(\Sigma \subset \lbrack -\beta ,\beta ])$ is a basis for $V_{\Lambda }$
whose elements are Delta ultrafunctions. Its dual basis $\left\{ \sigma
_{a}(x)\right\} _{a\in \Sigma }$ is called Sigma-basis. The set $\Sigma
\subset \lbrack -\beta ,\beta ]$ is called set of independent points.
\end{definition}

As we proved in \cite{belu2013}, Theorem 19, for every $q\in \lbrack -\beta
,\beta ]$ there exists a unique Delta ultrafunction centered in $q$. Let us
also note that, by saying that $\left\{ \sigma _{a}(x)\right\} _{a\in \Sigma
}$ is the dual basis of $\left\{ \mathbf{\delta }_{a}(x)\right\} _{a\in
\Sigma },$ we commit an abuse of language: in fact, in general, given a
basis $\left\{ e_{j}\right\} _{j=1}^{n}$ in a finite dimensional vector
space $V,$ the dual basis of $\left\{ e_{j}\right\} _{j=1}^{n}$ is the basis 
$\left\{ e_{j}^{\prime }\right\} _{j=1}^{n}$ of the dual space $V^{\prime }$
defined, for every $1\leq j,k\leq n$, by the following relation:%
\begin{equation*}
e_{j}^{\prime }\left[ e_{k}\right] =\delta _{jk}.
\end{equation*}%
When $V$ has a scalar product$\ \left( \cdot \ |\ \cdot \right) $ there
exists a base $g_{1},...,g_{n}$ of the space $V$ such that, for every $1\leq
j,k\leq n,$ we have%
\begin{equation*}
\left( g_{j}\ |\ e_{k}\right) =\delta _{jk}.
\end{equation*}

So $\left\{ e_{j}^{\prime }\right\} _{j=1}^{n}$ and $\{g_{j}\}_{j=1}^{n}$
can be identified, and $\{g_{j}\}_{j=1}^{n}$ will be called the dual basis
of $\left\{ e_{j}\right\} _{j=1}^{n}.$

In our case the scalar product that we consider is the extension of the $%
L^{2}$ scalar product to $V_{\Lambda }$, namely the scalar product such
that, for every $u,v\in V_{\Lambda },$ we have

\begin{equation*}
(u,v)=\int_{-\beta }^{\beta }u(x)v(x)dx.
\end{equation*}%
So a Sigma-basis is characterized by the fact that, $\forall a,b\in \Sigma ,$%
\begin{equation}
\int_{-\beta }^{\beta }\delta _{a}(x)\sigma _{b}(x)dx=\delta _{ab}.
\label{mimma}
\end{equation}

The existence of a Delta-basis (and, consequently, of a Sigma-basis) is an
immediate consequence of the following fact:

\begin{remark}
The set $\left\{ \mathbf{\delta }_{a}(x)|a\in \lbrack -\beta ,\beta
]\right\} $ generates all $V_{\Lambda }.$ In fact, let $G$ be the vector
space generated by the set $\left\{ \mathbf{\delta }_{a}(x)\ |\ a\in \lbrack
-\beta ,\beta ]\right\} $ and let us suppose that $G$ is properly included
in $V_{\Lambda }.$ Then the orthogonal $G^{\perp }$ of $G$ in $V_{\Lambda }$
contains a function $f\neq 0.$ But, since $f\in $ $G^{\perp },$ for every $%
a\in \lbrack -\beta ,\beta ]$ we have 
\begin{equation*}
f(a)=\int_{-\beta }^{\beta }f(x)\mathbf{\delta }_{a}(x)dx=0,
\end{equation*}%
so $f_{\upharpoonleft _{\lbrack -\beta ,\beta ]}}=0$ and this is absurd.
Thus the set $\left\{ \mathbf{\delta }_{a}(x)\ |\ a\in \lbrack -\beta ,\beta
]\right\} $ generates $V_{\Lambda },$ hence it contains a basis.
\end{remark}

Finally, let us recall the properties of a Sigma basis that we will use (see 
\cite{belu2013}, Theorem 22 for a proof of these results):

\begin{theorem}
\label{ytr}A Sigma-basis $\left\{ \sigma _{q}(x)\right\} _{q\in \Sigma }$
satisfies the following properties:

\begin{enumerate}
\item if $u\in V_{\Lambda }$ then%
\begin{equation*}
u(x)=\sum_{q\in \Sigma }u(q)\sigma _{q}(x);
\end{equation*}

\item if two ultrafunctions $u$ and $v$ coincide on a set of independent
points then they are equal;

\item if $\Sigma $ is a set of independent points and $a,b\in \Sigma $ then $%
\sigma _{a}(b)=\delta _{ab}$.
\end{enumerate}
\end{theorem}

For our aims, we need to fix an independent set $\Sigma $ that extends $%
\mathbb{R}\cup \{-\beta ,\beta \}$. This is possible, as the following
Theorem shows:

\begin{theorem}
There exists an independent set $\Sigma \subseteq \lbrack -\beta ,\beta ]$
such that 
\begin{equation*}
\mathbb{R}\cup \{-\beta ,\beta \}\subseteq \Sigma .
\end{equation*}
\end{theorem}

\textbf{Proof.\ }Given $a\in \mathbb{R}$ let%
\begin{equation*}
\Sigma _{a}=\{\Sigma \subseteq \lbrack -\beta ,\beta ]\mid \Sigma \text{ is
an independent set and }a,-\beta ,\beta \in \Sigma \}.
\end{equation*}

Each set $\Sigma _{a}$ is internal so, if we prove that the family $\{\Sigma
_{a}\}_{a\in \mathbb{R}}$ has the finite intersection property, we can
conclude by $\mathfrak{c}^{+}$-saturation (which holds, since we have chosen
to work in a $\left( 2^{\mathfrak{c}}\right) ^{+}-$saturated model).

Let $a_{1},...,a_{n}$ be distinct real numbers. To prove that $\Sigma
_{a_{1}}\cap ...\cap \Sigma _{a_{n}}\neq \emptyset $ it is sufficient to
show that the functions $\mathbf{\delta }_{a_{1}},...,\mathbf{\delta }%
_{a_{n}},\mathbf{\delta }_{-\beta },\mathbf{\delta }_{\beta }$ are linearly
independent (by duality, this fact entails that $\sigma _{a_{1}},...,\sigma
_{a_{n}},\sigma _{-\beta },\sigma _{\beta }$ are linearly independent, and
hence we have our thesis). We want to prove this fact.

First of all, $\mathbf{\delta }_{-\beta }$ and $\mathbf{\delta }_{\beta }$
are linearly independent, otherwise we would find an hyperreal number $\xi $
such that $\mathbf{\delta }_{\beta }=\xi \mathbf{\delta }_{-\beta },$ so $%
u(\beta )=\xi u(-\beta )$ for every ultrafunction $u,$ and this is clearly
false. For the general case let us suppose, by contrast, that 
\begin{equation*}
\mathbf{\delta }_{a_{1}}(x)=\sum_{i=2}^{n}c_{i}\mathbf{\delta }%
_{a_{i}}(x)+d_{1}\mathbf{\delta }_{-\beta }(x)+d_{2}\mathbf{\delta }_{\beta
}(x).
\end{equation*}%
Let $f\in \mathcal{C}_{0}^{1}(\mathbb{R})$ be such that $f(a_{1})\neq 0$
while $f(a_{i})=0$ for every $i=2,...,n$. Since $f\in \mathcal{C}_{0}^{1}(%
\mathbb{R})$ we have 
\begin{equation*}
\int_{-\beta }^{\beta }f^{\ast }(x)\mathbf{\delta }_{\beta
}(x)dx=\int_{-\beta }^{\beta }f^{\ast }(x)\mathbf{\delta }_{-\beta }(x)dx=0.
\end{equation*}%
Then%
\begin{eqnarray*}
0 &\neq &f^{\ast }(a_{1})=\int_{-\beta }^{\beta }f^{\ast }(x)\mathbf{\delta }%
_{a_{1}}(x)= \\
\int_{-\beta }^{\beta }f^{\ast }(x)\sum_{i=2}^{n}c_{i}\mathbf{\delta }%
_{a_{i}}(x)dx &=&\sum_{i=2}^{n}c_{i}\int_{-\beta }^{\beta }f^{\ast }(x)%
\mathbf{\delta }_{a_{i}}(x)dx=0,
\end{eqnarray*}%
which is absurd. $\square $

\bigskip

In the next section we will use an indipendent set of point $\Sigma $ to
define the notion of restricted ultrafunction. The algebra that we are
searching for will be precisely an algebra of restricted ultrafunctions.

\subsection{The Algebra of Restricted Ultrafunctions\label{dimostraz}}

Let us fix an independent set of points $\Sigma $ with $\mathbb{R}\cup
\{-\beta ,\beta \}\subseteq \Sigma .$ By point (1) in Proposition \ref{ytr}
it follows that every ultrafunction $u$ depends only on the values it
attains on an independent set of points; therefore, if $\mathfrak{I}\left(
\Sigma ,\mathbb{R}^{\ast }\right) $ is the family of internal functions $%
u:\Sigma \rightarrow \mathbb{R}^{\ast }$, then the operator of restriction $%
\Psi :V_{\Lambda }\rightarrow \mathfrak{I}\left( \Sigma ,\mathbb{R}^{\ast
}\right) $ given by 
\begin{equation*}
\Psi \left[ f\right] :=f_{\upharpoonleft _{\Sigma }}
\end{equation*}

is an isomorphism. The set $\mathfrak{I}\left( \Sigma ,\mathbb{R}^{\ast
}\right) $ will be denoted by $V(\Sigma ).$

\begin{definition}
The elements of $V(\Sigma )$ will be called restricted ultrafunctions.
\end{definition}

In order to simplify the notation, if $u$ is a restricted ultrafunction we
will write 
\begin{equation*}
\widetilde{u}:=\Psi ^{-1}\left[ u\right] .
\end{equation*}

Namely, if $\{\sigma _{a}(x)\}_{a\in \Sigma }$ is the Sigma-basis of $%
V_{\Sigma }$ associated to the independent set of points $\Sigma $, then%
\begin{equation*}
\widetilde{u}=\sum_{a\in \Sigma }u(a)\sigma _{a}(x).
\end{equation*}

The restricted ultrafunctions present the advantage that they form an
algebra with respect to the pointwise sum and product:%
\begin{equation*}
\left( f+g\right) (x)=f(x)+g(x);\ \left( f\cdot g\right) (x)=f(x)\cdot g(x).
\end{equation*}%
Moreover every restricted ultrafunction can be written as follows%
\begin{equation*}
u(x)=\sum_{a\in \Sigma }u(a)\delta _{ax},
\end{equation*}%
where $\delta _{ax}:\Sigma \rightarrow \{0,1\}$ is the usual Kronecker delta.

The spaces $V_{\Lambda }$ and $V(\Sigma )$ are isomorphic with respect to
many operations (e.g., with respect to the operations of sum and
multiplication by a constant) but not to all. This can be seen if we observe
that, when endowed with the pointwise multiplication, $V(\Sigma )$ is an
algebra while $V_{\Lambda }$ is not. In particular, if $u$ and $v$ are
restricted ultrafunctions, $\widetilde{u}\cdot \widetilde{v}$ is not in
general an extended ultrafunction, namely $\widetilde{u}\cdot \widetilde{v}%
\notin V_{\Lambda }$ and 
\begin{equation*}
\widetilde{u}\cdot \widetilde{v}\neq \widetilde{u\cdot v}.
\end{equation*}%
In any case, $\widetilde{u}\cdot \widetilde{v}$ and $\widetilde{u\cdot v}$
coincide on the points of $\Sigma $.

A nice feature of $V(\Sigma )$ is that it contains an extension of every
function $f\in \mathfrak{F}(\mathbb{R}):$

\begin{definition}
Given a function $f\in \mathfrak{F}(\mathbb{R}),$ its hyperfinite extension
(denoted by $f%
{{}^\circ}%
)$ is the restricted ultrafunction%
\begin{equation*}
f%
{{}^\circ}%
(x)=\sum_{a\in \Sigma }f^{\ast }(a)\delta _{ax}.
\end{equation*}
\end{definition}

We observe that, by definition, given any function $f\in \mathfrak{F}(%
\mathbb{R})$ we have%
\begin{equation*}
\widetilde{f%
{{}^\circ}%
}(x)=\sum_{a\in \Sigma }f^{\ast }(a)\sigma _{a}(x).
\end{equation*}%
So, in general, $\widetilde{f%
{{}^\circ}%
}(x)\neq f^{\ast }(x)$, even if for every $f\in \mathcal{C}^{1}(\mathbb{R})$
we have $\widetilde{f%
{{}^\circ}%
}(x)=f^{\ast }(x)$ (equivalently, for every $f\in \mathcal{C}^{1}(\mathbb{R}%
) $ we have $f%
{{}^\circ}%
=\Psi (f^{\ast })$).

We now introduce a scalar product on $V(\Sigma )$ that will play a central
role in what follows:

\begin{definition}
We denote by $\left\langle \cdot ,\cdot \right\rangle :V(\Sigma )\rightarrow 
\mathbb{R}^{\ast }$ the scalar product such that, for every $u,v\in V(\Sigma
),$ we have%
\begin{equation*}
\left\langle u,v\right\rangle =\int_{-\beta }^{\beta }\widetilde{u}(x)\cdot 
\widetilde{v}(x)dx.
\end{equation*}
\end{definition}

Notice that, in general, $\left\langle u,v\right\rangle \neq \int_{-\beta
}^{\beta }\widetilde{u\cdot v}(x)\ dx;$ in fact%
\begin{equation*}
\int_{-\beta }^{\beta }\widetilde{uv}(x)\ dx=\sum_{a\in \Sigma }u(a)v(a)\eta
_{a},
\end{equation*}%
while%
\begin{equation*}
\left\langle u,v\right\rangle =\sum_{a,b\in \Sigma }u(a)v(b)\eta _{ab},
\end{equation*}%
where, for every $a,b\in \Sigma $, we set%
\begin{equation*}
\eta _{a}=\int_{-\beta }^{\beta }\sigma _{a}(x)dx;\text{ \ \ \ }\eta
_{ab}=\int_{-\beta }^{\beta }\sigma _{a}(x)\sigma _{b}(x)dx.
\end{equation*}

Nevertheless, given any $f,g\in \mathcal{C}^{1}(\mathbb{R}),$ we have%
\begin{equation*}
\left\langle f%
{{}^\circ}%
,g%
{{}^\circ}%
\right\rangle =\int_{-\beta }^{+\beta }f^{\ast }(x)g^{\ast }(x)dx
\end{equation*}%
so, in particular, if $f,g\in \mathcal{C}_{0}^{1}(\mathbb{R})$ then%
\begin{equation*}
\left\langle f%
{{}^\circ}%
,g%
{{}^\circ}%
\right\rangle =\int f(x)g(x)dx.
\end{equation*}%
We use this scalar product to define the derivative \textrm{D}$:V(\Sigma
)\rightarrow V(\Sigma )$ by duality:

\begin{definition}
The derivative of a restricted ultrafunction $u$ (denoted by $\mathrm{D}u)$
is the unique restricted ultrafunction such that, $\forall \varphi \in
V(\Sigma ),$ we have 
\begin{equation*}
\left\langle \mathrm{D}u,\varphi \right\rangle =\int_{-\beta }^{\beta
}\partial ^{\ast }\widetilde{u}(x)\widetilde{\varphi }(x)dx.
\end{equation*}
\end{definition}

Let us observe that, since $\int_{-\beta }^{\beta }\partial ^{\ast }%
\widetilde{u}(x)\widetilde{\varphi }(x)dx=\int_{-\beta }^{\beta }D\widetilde{%
u}(x)\widetilde{\varphi }(x)dx,$ then 
\begin{equation*}
\widetilde{\mathrm{D}u}=D\widetilde{u}.
\end{equation*}%
So we can equivalently define \textrm{D} as follows: 
\begin{equation*}
\mathrm{D}=\Psi \circ D\circ \Psi ^{-1}=\Psi \circ P_{V_{\Lambda }}\circ
\partial ^{\ast }\circ \Psi ^{-1}.
\end{equation*}%
In particular \textrm{D}$f%
{{}^\circ}%
=(\partial f)%
{{}^\circ}%
$ whenever $f\in \mathcal{C}^{2}(\mathbb{R}).$

By combining Theorem \ref{facundo} with the definitions of the scalar
product $\left\langle \cdot ,\cdot \right\rangle $ and of the operator 
\textrm{D} we obtain the following result:

\begin{theorem}
\label{wlr}For every $u,v\in V(\Sigma )$ we have 
\begin{equation*}
\left\langle \mathrm{D}u(x),v(x)\right\rangle =-\left\langle u(x),\mathrm{D}%
v(x)\right\rangle +[uv]_{-\beta }^{\beta }.
\end{equation*}
\end{theorem}

\textbf{Proof. }Let us compute $\left\langle \mathrm{D}u(x),v(x)\right%
\rangle :$%
\begin{eqnarray*}
\left\langle \mathrm{D}u(x),v(x)\right\rangle &=&\int_{-\beta }^{\beta }%
\widetilde{\mathrm{D}u(x)}\widetilde{v(x)}dx=\int_{-\beta }^{\beta }D%
\widetilde{u(x)}\widetilde{v(x)}dx= \\
-\int_{-\beta }^{\beta }\widetilde{u(x)}D\widetilde{v(x)}dx+[\widetilde{u}%
\widetilde{v}]_{-\beta }^{\beta } &=&-\int_{-\beta }^{\beta }\widetilde{u(x)}%
\widetilde{\mathrm{D}v(x)}dx+[\widetilde{u}\widetilde{v}]_{-\beta }^{\beta }
\\
&=&-\left\langle u(x),\mathrm{D}v(x)\right\rangle +[uv]_{-\beta }^{\beta }.%
\text{ \ \ \ \ \ \ \ \ \ \ \ \ \ \ \ \ \ \ }\square
\end{eqnarray*}

\bigskip

Now we want to define a (in some sense canonical) embedding of distributions 
\begin{equation*}
\Phi :\mathcal{D}^{\prime }(\mathbb{R})\rightarrow V(\Sigma ).
\end{equation*}%
A known representation theorem for distributions states that for every
distribution $T\in \mathcal{D}^{\prime }(\mathbb{R})$ and for every compact
set $[a,b]\subseteq \mathbb{R}$ there exist $f\in \mathcal{C}^{1}\left( 
\mathbb{R}\right) $ and $k\in \mathbb{N}$ such that $T_{\upharpoonleft
_{\lbrack a,b]}}=\partial ^{k}f$. By transfer we deduce that, in particular,
for every distribution $T\in \mathcal{D}^{\prime }(\mathbb{R})$ there exists 
$\varphi \in \mathcal{C}^{1}(\mathbb{R)}^{\ast }$ and $k\in \mathbb{N}^{\ast
}$ such that $T_{\upharpoonleft _{\lbrack -\beta ,\beta ]}}^{\ast }=\partial
^{k}\varphi .$ So the following definition makes sense:

\begin{definition}
Given $T\in \mathcal{D}^{\prime }(\mathbb{R})$ we define%
\begin{equation*}
d_{T}=\min \{k\in \mathbb{N}^{\ast }\mid \exists \varphi \in \mathcal{C}^{1}(%
\mathbb{R)}^{\ast }\text{such that }T_{\upharpoonleft _{\lbrack -\beta
,\beta ]}}^{\ast }=\partial ^{k}\varphi \}
\end{equation*}%
and 
\begin{equation*}
R_{T}=\{\varphi \in \mathcal{C}^{1}(\mathbb{R)}^{\ast }|T_{\upharpoonleft
_{\lbrack -\beta ,\beta ]}}^{\ast }=\partial ^{d_{T}}\varphi \}.
\end{equation*}
\end{definition}

It is known that if the weak derivative of a continuous function is zero on
an interval $[a,b]$ then the continuous function is constant on $[a,b]$. It
is easy to generalize this result and to prove the following Lemma:

\begin{lemma}
\label{javierzanetti}Let $f\in \mathcal{C}^{0}(\mathbb{R)}$, $k\in \mathbb{N}
$. If $\partial ^{k}f_{\upharpoonleft _{\lbrack a,b]}}=0$ then there exists
a polynomial $P(x)$, with $\deg (P(x))<k,$ such that $f(x)_{\upharpoonleft
_{\lbrack a,b]}}=P(x).$
\end{lemma}

This standard result allows us to prove the following

\begin{lemma}
\label{reginaldo}If $\varphi _{1},\varphi _{2}\in R_{T}$ then then there
exists a polynomial $P(x)\in \mathcal{C}^{1}(\mathbb{R)}^{\ast }$, with $%
\deg (P(x))<d_{T},$ such that $\left( \varphi _{1}-\varphi _{2}\right)
_{\upharpoonleft _{\lbrack -\beta ,\beta ]}}=P(x).$
\end{lemma}

\textbf{Proof. }By construction, $\varphi _{1}-\varphi _{2}\in \mathcal{C}%
^{1}(\mathbb{R)}^{\ast }$ and the $d_{T}$-th weak derivative of $\varphi
_{1}-\varphi _{2}$ is zero on $[-\beta ,\beta ]$. We apply the nonstandard
version of Lemma \ref{javierzanetti} obtained by transfer and we deduce the
thesis. \ $\square $

\bigskip

Let us also observe that $d_{\partial T}=d_{T}+1$ and that by Lemma \ref%
{reginaldo} it follows that 
\begin{equation*}
R_{\partial T}=\{\varphi _{T}+rx^{d_{T}}\mid \varphi _{T}\in R_{T},r\in 
\mathbb{R}^{\ast }\}.
\end{equation*}

\begin{theorem}
\label{birra}There exists a hyperfinite set $H$ such that $R_{T}\subseteq
Span(H)$ for every $T\in \mathcal{D}^{\prime }(\mathbb{R}).$
\end{theorem}

\textbf{Proof.}\ By saturation, the intersection%
\begin{equation*}
\bigcap_{T\in D^{\prime }}[d_{T},+\infty )
\end{equation*}

is nonempty. Let $\alpha \in \bigcap\limits_{T\in D^{\prime }}[d_{T},+\infty
).$ For every $T\in \mathcal{D}^{\prime }(\mathbb{R})$ let $\varphi _{T}\in
R_{T}\mathcal{\ }$and let%
\begin{equation*}
H_{T}=\{H\subseteq \mathcal{C}^{1}(\mathbb{R)}^{\ast }\mid H\text{ is
hyperfinite and }1,x,...,x^{\alpha },\varphi _{T}\in H\}.
\end{equation*}

For every $T\in \mathcal{D}^{\prime }(\mathbb{R})$ the set $H_{T}$ is
nonempty; moreover, the family $\{H_{T}\}_{T\in D^{\prime }}$ has the finite
intersection property since, for every $H_{1}\in H_{T_{1}},H_{2}\in
H_{T_{2}},$ we have $H_{1}\cup H_{2}\in H_{T_{1}}\cap H_{T_{2}}.$ Then by
saturation we have 
\begin{equation*}
\bigcap_{T\in D^{\prime }}H_{T}\neq \emptyset .
\end{equation*}%
Let $H\in \bigcap\limits_{T\in D^{\prime }}H_{T}.$ For every $T\in \mathcal{D%
}^{\prime }(\mathbb{R})$ we have that $1,x,...,x^{d_{T}},\varphi _{T}\in H$
hence, by Lemma \ref{reginaldo}, we conclude that $R_{T}\in Span\left(
H\right) .$ $\square $

\bigskip

Now let $H$ and $\alpha $ be given as in Theorem \ref{birra}. Let%
\begin{equation*}
\widetilde{H}=\{P^{s}(h)\mid 0\leq s\leq \alpha \text{ and }h\in H\},
\end{equation*}

where $P$ denotes the operator that maps a function in $\mathcal{C}^{1}(%
\mathbb{R)}^{\ast }$ to (one of) its primitive with respect to $\partial .$
From now on we consider $\Lambda \subseteq \mathbb{C}^{1}(\mathbb{R)}^{\ast
} $ to be a hyperfinite set with 
\begin{equation*}
\widetilde{H}\cup \mathfrak{F}(\mathbb{R},\mathbb{R})^{\sigma }\subseteq
\Lambda
\end{equation*}%
and we construct our model by mean of this hyperfinite set $\Lambda .$ Let
us note that, as a consequence of this choice, we have that 
\begin{equation*}
\forall T\in \mathcal{D}^{\prime }(\mathbb{R})\text{ }\forall \varphi
_{T}\in R_{T}\text{ }\forall s\in \mathbb{N}^{\ast }\cap \lbrack 0,\alpha ]%
\text{ }P^{s}(\varphi _{T})\text{ is an ultrafunction.}
\end{equation*}%
In particular every polynomial $P(x)$ with $\deg (P(x))\leq \alpha $ is an
ultrafunction.

\begin{lemma}
\label{crescenza}Let $P(x)\in V_{\Lambda }$ be a polynomial and let $\deg
(P(x))<\alpha$. Then $\mathrm{D}^{k+1}(\Psi (P(x)))=0.$
\end{lemma}

\textbf{Proof.} Since $\deg (P)<\alpha ,$ $P(x),\partial P(x),...,\partial
^{k+1}P(x)$ are ultrafunctions, so we deduce that $D^{i}P(x)=\partial
^{i}P(x)$ for every $0\leq i\leq k+1.$ Then $D^{k+1}P(x)=\partial
^{k+1}P(x)=0,$ and we obtain the thesis by recalling that $\mathrm{D}%
^{k+1}(\Psi (P(x)))=\Psi (D^{k+1}P(x))$. \ $\square $

\begin{definition}
We denote by $\Phi :\mathcal{D}^{\prime }\left( \mathbb{R}\right)
\rightarrow V(\Sigma )$ the function such that for every $T\in \mathcal{D}%
^{\prime }\left( \mathbb{R}\right) $%
\begin{equation}
\Phi (T)=\mathrm{D}^{d_{T}}(\Psi (\varphi _{T})),
\end{equation}%
where $\varphi _{T}\in R_{T}.$
\end{definition}

Lemma \ref{crescenza} entails that $\Phi $ is well defined since it does not
depend on the particular choice of $\varphi_{T}\in R_{T}$. The function $%
\Phi $ has a few important properties:

\begin{theorem}
\label{dualita} We have the following properties:

\begin{enumerate}
\item if $f\in \mathcal{C}^{1}(\mathbb{R})$ then $\Phi (T_{f})=f%
{{}^\circ}%
;$

\item $\forall T\in \mathcal{D}^{\prime }(\mathbb{R}),\ \forall \varphi \in 
\mathcal{D}(\mathbb{R}), T(\varphi )=\left\langle \Phi (T),\varphi 
{{}^\circ}%
\right\rangle ;$

\item the following diagram commutes: 
\begin{equation*}
\begin{array}{ccc}
\mathcal{D}^{\prime }(\mathbb{R}) & \overset{\partial }{\longrightarrow } & 
\mathcal{D}^{\prime }(\mathbb{R}) \\ 
\downarrow \Phi &  & \downarrow \Phi \\ 
V(\Sigma ) & \overset{\mathrm{D}}{\longrightarrow } & V(\Sigma )%
\end{array}%
\end{equation*}%
where $\partial $ is the usual distributional derivative;

\item the restriction of $\cdot $ to $\mathcal{C}^{1}(\mathbb{R})$ agrees
with the pointwise product, namely if $f,g\in \mathcal{C}^{1}(\mathbb{R})\ $%
then%
\begin{equation*}
\Phi (T_{fg})=\Phi (T_{f})\cdot \Phi \left( T_{g}\right) .
\end{equation*}
\end{enumerate}
\end{theorem}

\textbf{Proof.} 1) If $f\in \mathcal{C}^{1}(\mathbb{R})$ then $T_{f}=f$. So $%
d_{f}=0$ and $R_{T}=\{f^{\ast }\}$ hence, by definition, $\Phi (T_{f})=\Psi
(f^{\ast })=f%
{{}^\circ}%
.$

2) Let $T_{\upharpoonleft _{\lbrack -\beta ,\beta ]}}^{\ast }=\partial
^{d_{T}}f.$ We compute $\left\langle \Phi (T),\varphi 
{{}^\circ}%
\right\rangle :$%
\begin{eqnarray*}
\left\langle \Phi (T),\varphi 
{{}^\circ}%
\right\rangle &=&\int_{-\beta }^{\beta }\widetilde{\Phi (T)}\cdot \widetilde{%
\varphi 
{{}^\circ}%
}dx= \\
\int_{-\beta }^{\beta }\widetilde{\mathrm{D}^{d_{T}}(\Psi (f))}\cdot \varphi
^{\ast }dx &=&\int_{-\beta }^{\beta }D^{d_{T}}f\cdot \varphi ^{\ast }dx.
\end{eqnarray*}%
Now by Proposition \ref{derivatona} it follows that%
\begin{equation*}
\int_{-\beta }^{\beta }D^{d_{T}}f\cdot \varphi ^{\ast
}dx=(-1)^{d_{T}}\int_{-\beta }^{\beta }f\cdot \partial ^{d_{T}}\varphi
^{\ast }dx.
\end{equation*}%
So%
\begin{eqnarray*}
\left\langle \Phi (T),\varphi 
{{}^\circ}%
\right\rangle &=&\int_{-\beta }^{\beta }D^{d_{T}}f\cdot \varphi ^{\ast }dx=
\\
(-1)^{d_{T}}\int_{-\beta }^{\beta }f\cdot \partial ^{d_{T}}\varphi ^{\ast
}dx&=&T^{\ast }[\varphi ^{\ast }]= \\
\left( T[\varphi ]\right) ^{\ast }=T[\varphi].
\end{eqnarray*}

3) Let $T\in \mathcal{D}^{\prime }(\mathbb{R}),$ $T_{\upharpoonleft
_{\lbrack -\beta ,\beta ]}}^{\ast }=\partial ^{d_{T}}f,$ $f\in R_{T}$. Let
us compute $\mathrm{D}(\Phi (T)):$%
\begin{eqnarray*}
\mathrm{D}(\Phi (T)) &=&\mathrm{D}\left( \Psi (D^{d_{T}}(f))\right) = \\
\Psi (D^{d_{T}+1}(f)) &=&\Phi (\partial T),
\end{eqnarray*}

since $d_{\partial T}=d_{T}+1$ and $f\in R_{\partial T}.$

4) Since $f,g\in \mathcal{C}^{1}(\mathbb{R})$ then $\Phi (T_{fg})=(fg)%
{{}^\circ}%
=f%
{{}^\circ}%
g%
{{}^\circ}%
=\Phi (T_{f})\cdot \Phi (T_{g}).$ $\square $

\begin{corollary}
\label{shotgun}$\Phi :\mathcal{D}^{\prime }\left( \mathbb{R}\right)
\rightarrow V(\Sigma )$ is an embedding of vector spaces.
\end{corollary}

\textbf{Proof.} $\Phi $ is injective: if $T_{1}\neq T_{2}$ are distributions
then there is a test function $\varphi \in \mathcal{D}(\mathbb{R})$ such
that $T_{1}[\varphi ]\neq T_{2}[\varphi ].$ In particular, 
\begin{equation*}
\left\langle \Phi (T_{1}),\varphi 
{{}^\circ}%
\right\rangle =T_{1}[\varphi ]\neq T_{2}[\varphi ]=\left\langle \Phi
(T_{2}),\varphi 
{{}^\circ}%
\right\rangle ,
\end{equation*}%
hence $\Phi (T_{1})\neq \Phi (T_{2}).$

$\Phi $ is a linear map: let $T_{1\upharpoonleft _{\lbrack -\beta ,\beta
]}}^{\ast }=\partial ^{d_{T_{1}}}f,$ $T_{2\upharpoonleft _{\lbrack -\beta
,\beta ]}}^{\ast }=\partial ^{d_{T_{2}}}g,$ $f,g$ ultrafunctions. Let us
suppose that $d_{T_{1}}=d_{T_{2}}+s.$ Necessarily, $0\leq s\leq \alpha $
since both $d_{T_{1}},d_{T_{2}}\leq \alpha .$ So%
\begin{equation*}
\left( T_{1}+T_{2}\right) _{\upharpoonleft _{\lbrack -\beta ,\beta ]}}^{\ast
}=\partial ^{d_{T_{1}}}(f+P^{s}(g)),
\end{equation*}%
therefore%
\begin{eqnarray*}
\Phi (T_{1}+T_{2}) &=&\mathrm{D}^{d_{T_{1}}}(\Psi (f+P^{s}(g)))= \\
&&\mathrm{D}^{d_{T_{1}}}(\Psi (f))+\mathrm{D}^{d_{T_{1}}}(\Psi (P^{s}(g))).
\end{eqnarray*}%
Now, by definition $\mathrm{D}^{d_{T_{1}}}(\Psi (f))=\Phi (T_{1}).$ Moreover%
\begin{eqnarray*}
\mathrm{D}^{d_{T_{1}}}(\Psi (P^{s}(g))) &=&\mathrm{D}^{d_{T_{2}}+s}(\Psi
(P^{s}(g)))= \\
\mathrm{D}^{d_{T_{2}}}(\mathrm{D}^{s}(\Psi (P^{s}(g)))) &=&\mathrm{D}%
^{d_{T_{2}}}(\Psi (D^{s}(P^{s}(g)))).
\end{eqnarray*}%
By our choice of $\Lambda $ in the construction of the space of
ultrafunctions we have that $g,P(g),...,P^{s}(g)$ are ultrafunctions. So $%
D^{s}(P^{s}(g))=\partial ^{s}(P^{s}(g))=g,$ hence%
\begin{eqnarray*}
\mathrm{D}^{d_{T_{1}}}(\Psi (P^{s}(g))) &=&\mathrm{D}^{d_{T_{2}}}(\Psi
(D^{s}(P^{s}(g))))= \\
\mathrm{D}^{d_{T_{2}}}(\Psi (g)) &=&\Phi (T_{2}).
\end{eqnarray*}

To prove that 
\begin{equation*}
\Phi (rT)=r\Phi (T)
\end{equation*}%
it is sufficient to observe that $d_{rT}=d_{T}$ and $R_{rT}=rR_{T}.$ This
proves that $\Phi $ is linear. In particular if $r=0$ we get that the image
of the zero distribution is the zero ultrafunction, as expected. $\square $%
\bigskip

We are now ready to prove Theorem \ref{gonerilla}:

\bigskip

\textbf{Proof of Theorem 1: }Let us pose $\mathfrak{A}=V(\Sigma ),$ and let
us consider $\Phi ,$\QTR{cal}{\textrm{D}}$,\left\langle \cdot ,\cdot
\right\rangle $ as introduced in this section. Then ($\mathfrak{A}$-0)
follows by the definition of $\mathfrak{A};$ ($\mathfrak{A}$-1), ($\mathfrak{%
A}$-2) and ($\mathfrak{A}$-3) have been proved in Theorem \ref{dualita} and
Corollary \ref{shotgun}; ($\mathfrak{A}$-4) has been proved in Theorem \ref%
{wlr} and ($\mathfrak{A}$-5) follows immediatly by the definition of $f%
{{}^\circ}%
$. $\square $

\bigskip

\textbf{Acknowledgement: }The authors would like to thank the referee for
his careful reading of the paper.

\end{document}